\documentclass{article}
\usepackage{amsmath,amssymb,amsthm}

\newtheorem{theorem}{Theorem}
\newtheorem{lemma}{Lemma}

\newtheorem{remark}{Remark}

\numberwithin{theorem}{section}
\numberwithin{definition}{section}
\numberwithin{lemma}{section}
\numberwithin{proposition}{section}
\numberwithin{remark}{section}

\begin{document}
\title{Boolean representation of manifolds and functions}
	\author{Sergei~Ovchinnikov \\
	Mathematics Department\\
	San Francisco State University\\
	San Francisco, CA 94132\\
	sergei@sfsu.edu} 
\date{\today} 
\maketitle


\begin{abstract}
\noindent
It is shown that a smooth $n$--dimensional manifold in $\mathbb{R}^n$ with a boundary admits a Boolean representation in terms of closed half spaces defined by the tangent hyperplanes at the points on its boundary. A similar result is established for smooth functions on closed convex domains in $\mathbb{R}^n$. The latter result is considered as a form of the Legendre transform.
\end{abstract}

\section{Introduction}
Let $U$ be a set, $X$ be a subset of $U$, and $\mathcal{X}=\{X_i\}_{i\in I}$ be a family of subsets of $U$. We say that $X$ admits a \emph{Boolean representation in terms of $\mathcal{X}$} if
\begin{equation} \label{BooleanDomain}
X = \bigcup_{j\in J}\bigcap_{i\in S_j} X_i
\end{equation}
for some family $\{S_j\}_{j\in J}$ of subsets of $I$.

Similarly, let $f$ be a real valued function on a set $\Omega$ and $\mathcal{G}=\{g_i\}_{i\in I}$ be a family of real valued functions on $\Omega$. We say that $f$ admits a \emph{Boolean representation in terms of $\mathcal{G}$} if
\begin{equation} \label{BooleanFunction}
f(x) = \sup_{j\in J}\,\inf_{i\in S_j} g_i(x),\quad\forall x\in\Omega.
\end{equation}
for some family $\{S_j\}_{j\in J}$ of subsets of $I$.

The goal of this paper is to establish Boolean representations for smooth domains in $\mathbb{R}^n$ ($n\geq 2$) and smooth real valued functions on closed convex subsets of $\mathbb{R}^n$. The paper is organized as follows.

In Section~2 we show that a smooth $n$--dimensional manifold in $\mathbb{R}^n$ with a boundary admits a Boolean representation (\ref{BooleanDomain}) in terms of a family of closed half subspaces of $\mathbb{R}^n$.

In Section~3 we establish a Boolean representation (\ref{BooleanFunction}) for smooth real valued functions on closed convex domains in $\mathbb{R}^n$. Although this result can be established using the method developed in Section~2, we prefer to present an alternative proof that gives an explicit formula for the family $\mathcal{G}$ in (\ref{BooleanFunction}).

Finally, in Section~4, we investigate the relation between the Boolean representation established in Section~3 and the Legendre transform.

\section{Smooth manifolds}

Let $\Gamma\subset\mathbb{R}^n$ be a smooth $n$--dimensional manifold with a boundary, i.e., each $x\in\Gamma$ has a neighborhood diffeomorphic to an open subset of a closed half space in $\mathbb{R}^n$. For each point $x$ in the boundary $\partial\Gamma$ we denote by $Q_x$ the closed half space consisting of all tangent and `inward' vectors at $x$.

Let $a$ be a point in $\Gamma$ and $R$ be a closed ray with the origin at $a$. The connected component of $\Gamma\cap R$ in $R$ is a closed interval $[a,a_R]$ in $R$ (it is possible that $a_R=a$). In other words, $[a,a_R]$ is the set of all points in $\Gamma$ that are `visible' from $a$ in the `direction' $R$. Clearly, $a_R\in\partial\Gamma$ and $[a,a_R]\subset Q_{a_R}$.

Let $\mathcal{R}_a$ be the set of all closed rays with the origin at $a$. Then $\bigcap_{R\in\mathcal{R}_a}Q_{a_R}$ is a  closed convex subset of $\Gamma$ containing $a$. Thus we have
\begin{equation*}
\Gamma = \bigcup_{a\in\Gamma}\,\bigcap_{R\in\mathcal{R}_a}Q_{a_R}.
\end{equation*}

The above argument proves the following theorem.

\begin{theorem} \label{BooleanDomainTheorem}
Let $\Gamma\subset\mathbb{R}^n$ be a smooth $n$--dimensional manifold with the boundary $\partial\Gamma$ and let $\mathcal{Q}$ be the family of all closed half spaces $Q_x,\; x\in\partial\Gamma$. There is a family $\{\mathcal{Q}_i\}_{i\in I}$ of subsets of $\mathcal{Q}$ such that $\Gamma$ admits a Boolean representation
\begin{equation} \label{BDrepresentation}
\Gamma = \bigcup_{i\in I}\,\bigcap_{Q_x\in\mathcal{Q}_i}Q_x.
\end{equation}
\end{theorem}

\begin{remark}\normalfont
Suppose $\Gamma$ is convex. Then $\Gamma = \bigcap_{x\in\partial\Gamma}Q_x$ which is clearly a special case of (\ref{BDrepresentation}). If $\Gamma$ is concave, then $\Gamma = \bigcup_{x\in\partial\Gamma}Q_x$ which is again a special case of (\ref{BDrepresentation}).
\end{remark}

\begin{remark}\normalfont
Representations similar to (\ref{BDrepresentation}) are well--known in the area of constructive solid geometry (CSG). Namely, any simple polytope can be represented by a Boolean formula based on the half spaces supporting the faces of the polytope \cite{dD93}.
\end{remark}

\section{Smooth functions}
We use notation $x=(x_1,\ldots,x_n)$ for points in $\mathbb{R}^n$. A closed domain in $\mathbb{R}^n$ is the closure of an open set.

Let $f$ be a smooth function on a closed domain $\Omega$ in $\mathbb{R}^n$. For a point $t\in\Omega$ we define
\begin{equation} \label{tangent-space}
g_t(x) = \langle\nabla f(t),x-t\rangle+f(t),\quad x\in\Omega,
\end{equation}
where $\nabla f(t)$ is the gradient vector of $f$ at $t$ and $\langle\cdot,\cdot\rangle$ is the standard inner product in $\mathbb{R}^n$.
Geometrically, the graphs of these affine linear functions are tangent hyperplanes to the graph of $f$.

In this section we prove the following theorem.

\begin{theorem} \label{main-theorem}
Given closed convex domain $\Omega$ and $f\in C^1(\Omega)$, there exists a family $\{S_j\}_{j\in J}$ of subsets of $\Omega$ such that
\begin{equation} \label{Boolean-representation}
f(x) = \sup_{j\in J}\,\inf_{t\in S_j} g_t(x),\quad\forall x\in\Omega.
\end{equation}
\end{theorem}

First we prove two technical lemmas.

\begin{lemma} \label{Lemma1}
Let $h$ be a differentiable function on $[0,1]$. There exists $\lambda_0\in[0,1]$ such that
\begin{equation} \label{Eq1}
h'(\lambda_0)(-\lambda_0)+h(\lambda_0)\leq h(0)\quad\text{and}\quad h'(\lambda_0)(1-\lambda_0)+h(\lambda_0)\geq h(1).
\end{equation}
\end{lemma}

\begin{proof}
Let $m=h(1)-h(0)$. If $h'(0)\geq m$, then $\lambda_0=0$ satisfies both inequalities. Thus we may assume that $h'(0)<m$. Similarly, if $h'(1)\geq m$, then $\lambda_0=1$ satisfies both inequalities and we may assume that $h'(1)<m$.

Consider function $H(\lambda)=h(\lambda)-m\lambda-h(0)$. We have
\begin{equation*}
H(0)=H(1)=0\quad\text{and}\quad H'(0)<0,\;H'(1)<0.
\end{equation*}
It follows that $H$ is negative in some neighborhood of $0$ and positive in some neighborhood of $1$. Hence the set $U=\{\lambda\in(0,1):H(\lambda)=0\}$ is a nonempty closed subset of $(0,1)$. Let $\lambda_0=\inf U$. Then $H(\lambda_0)=0$ and $H'(\lambda_0)\geq 0$, i.e.,
\begin{equation*}
h(\lambda_0)=m\lambda_0+h(0)=m(\lambda_0-1)+h(1)
\end{equation*}
and $h'(\lambda_0)\geq m$. We have
\begin{equation*}
h'(\lambda_0)(-\lambda_0)+h(\lambda_0)=h'(\lambda_0)(-\lambda_0)+m\lambda_0+h(0) \leq h(0)
\end{equation*}
and
\begin{equation*}
h'(\lambda_0)(1-\lambda_0)+h(\lambda_0)=h'(\lambda_0)(1-\lambda_0)+m(\lambda_0-1)+h(1)\geq h(1).
\end{equation*}
\end{proof}

\begin{lemma} \label{Lemma2}
Let $f\in C^1(\Omega)$. For any given $a,b\in\Omega$ there exists $c\in\Omega$ such that
\begin{equation*}
g_c(a)\leq f(a)\quad\text{and}\quad g_c(b)\geq f(b).
\end{equation*}
\end{lemma}

\begin{proof}
Let $h(\lambda)=f((1-\lambda)a+\lambda b)$ for $\lambda\in[0,1]$. By Lemma~\ref{Lemma1}, there is $\lambda_0\in[0,1]$ satisfying inequalities~(\ref{Eq1}). Let $c=(1-\lambda_0)a+\lambda_0 b$. We have
\begin{align*}
g_c(a)&=\langle\nabla f(c),a-c\rangle+f(c) = (-\lambda_0)\langle\nabla f(c),b-a\rangle+f(c) \\
&=(-\lambda_0)h'(\lambda_0)+h(\lambda_0) \leq h(0) = f(a)
\end{align*}
and
\begin{align*}
g_c(b)&=\langle\nabla f(c),b-c\rangle+f(c) = (1-\lambda_0)\langle\nabla f(c),b-a\rangle+f(c) \\
&=(1-\lambda_0)h'(\lambda_0)+h(\lambda_0) \geq h(1) = f(b).
\end{align*}

\end{proof}

Now we proceed with the proof of Theorem~\ref{main-theorem}. For a given $u\in\Omega$, we define $S_u=\{t\in\Omega:g_t(u)\geq f(u)\}$ and 
\begin{equation*}
f_u(x)=\inf_{t\in S_u} g_t(x),\quad x\in\Omega.
\end{equation*}
which is well--defined, since $f\in C^1(\Omega)$. 

By Lemma~\ref{Lemma2}, for given $x,u\in\Omega,\;x\not= u$, there exists $v\in\Omega$ such that $g_v(u)\geq f(u)$ and $g_v(x)\leq f(x)$. Hence, $v\in S_u$ and
\begin{equation*}
f_u(x)=\inf_{t\in S_u} g_t(x) \leq g_v(x) \leq f(x).
\end{equation*}
In addition,
\begin{equation*}
f_u(u) = \inf_{t\in S_u} g_t(u) = g_u(u)=f(u),
\end{equation*}
since $u\in S_u$ and $g_t(u)\geq f(u)$ for $t\in S_u$. Therefore we have
\begin{equation*}
\sup_{u\in\Omega}\,\inf_{t\in S_u} g_t(x)=\sup_{u\in\Omega}f_u(x)=f(x),
\end{equation*}
which completes the proof.

\begin{remark}\normalfont
Let $f$ be a strictly convex function. Then
\begin{equation*}
f(x)= \sup_{u\in\Omega} g_u(x),\quad x\in\Omega,
\end{equation*}
since $S_u=\{u\}$ for all $u\in\Omega$. Similarly, for a strictly concave $f$, $S_u=\Omega$ and we have
\begin{equation*}
f(x)= \inf_{u\in\Omega} g_u(x),\quad x\in\Omega.
\end{equation*}
Both facts are well known properties of convex (concave) functions.
\end{remark}

\begin{remark}\normalfont
Convexity of $\Omega$ is an essential assumption in Theorem~\ref{main-theorem}. Consider, for instance, the domain $\Omega$ in $\mathbb{R}^2$ which is a union of three triangles defined by the sets of their vertices as follows (see Fig.~1.):
\begin{align*}
&\Delta_1=\{(-1,0),(-1,-1),(0,0)\},\quad\Delta_2=\{(0,0),(1,1),(1,0)\}, \\
& \text{and}\quad\Delta_3=\{(-1,0),(1,0),(0,-1)\}.
\end{align*}

\medskip
\begin{center}
\begin{picture}(150,150)
	\put(0,75){\vector(1,0){150}}\put(155,75){$x_1$}
	\put(75,60){\vector(0,1){90}}\put(75,155){$x_2$}
	\put(75,5){\line(0,1){40}}
	\put(30,75){\line(0,1){45}}
	\put(120,75){\line(0,1){45}}
	\put(30,120){\line(1,-1){45}}
	\put(75,75){\line(1,1){45}}
	\put(30,75){\line(1,-1){45}}
	\put(75,30){\line(1,1){45}}
	\put(35,90){$\Delta_1$}
	\put(102,90){$\Delta_2$}
	\put(70,50){$\Delta_3$}
\end{picture}

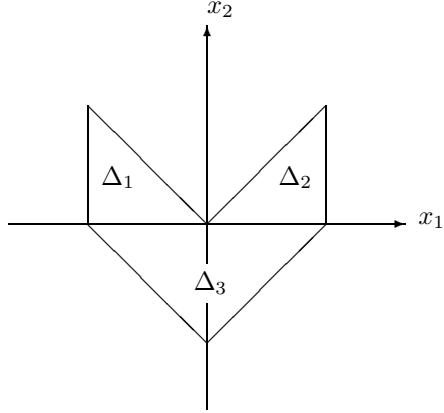
{\footnotesize Fig.~1}
\end{center}

\noindent
Let us define
\begin{equation*}
f(x) = \begin{cases}
	x_2^2, &\text{for $x\in\Delta_2,$} \\
	0, &\text{for $x\in\Delta_1\cup\Delta_3$.}
\end{cases}
\end{equation*}
Clearly, $f\in C^1(\Omega)$. Suppose $f$ has a Boolean representation
\begin{equation*}
f(x) = \sup_{j\in J}\,\inf_{t\in S_j} g_t(x),\quad\forall x\in\Omega,
\end{equation*}
for some family $\{S_j\}_{j\in J}$ of subsets of $\Omega$. Let $a$ be a point in the interior of $\Delta_2$. Since $f(a)>0$, there is $S_j$ such that $\inf \{g_t(a):t\in S_j\}>0$. Since $g_t(a)=0$ for $t\in\Delta_1\cup\Delta_3$, we have $S_j\subseteq\Delta_1$. Let $b=(-a_1,a_2)$. Then $g_t(b)=g_t(a)$, since $g_t(x)=2t_2 x_2-t_2^2$. Thus $\inf \{g_t(b):t\in S_j\}=\inf \{g_t(a):t\in S_j\}>0$ which contradicts $f(b)=0$.
\end{remark}

\begin{remark}\normalfont
A Boolean representation of piecewise linear functions in terms of their linear components similar to (\ref{Boolean-representation}) is obtained in \cite{sO00}.
\end{remark}

\section{The Legendre transform}

Let $f$ be a strictly convex smooth function on a compact convex domain $\Omega\subset\mathbb{R}^n$. Then (cf.~Remark~3.1)
\begin{equation} \label{convex-representation}
f(x) = \sup_{t\in\Omega} g_t(x) = \sup_{t\in\Omega}\,\{\langle\nabla f(t),x\rangle-[\langle\nabla f(t),t\rangle-f(t)]\},\quad\forall x\in\Omega.
\end{equation}
Let us introduce variables
\begin{align}
&p = \nabla f(t), \label{p} \\
&H=\langle\nabla f(t),t\rangle-f(t) \label{H}
\end{align}
Since $f$ is strictly convex, equation~(\ref{p}) defines a one--to--one mapping of $\Omega$ onto $\Omega'=\nabla f(\Omega)$ and we can express $t$ in terms of $p$ in~(\ref{convex-representation}) to obtain the following representation of $f(x)$ in terms of its \emph{Legendre transform} $H(p)$ (cf.~\cite{iGsF63}):
\begin{equation*}
f(x) = \sup_{p\in\Omega'}\{\langle p,x\rangle-H(p)\}
\end{equation*}

In general, let $\Gamma$ denote the graph of a smooth function $f:\Omega\rightarrow\mathbb{R}$. Then $\Gamma$ is the envelope of the set of its tangent hyperplanes. The \emph{Legendre transform}~\cite{lY69} of $\Gamma$ is the surface defined parametrically by~(\ref{p}) and~(\ref{H}) in the $(n+1)$--dimensional $(p,H)$--space. One can view our Boolean representation~(\ref{Boolean-representation}) as a representation of an arbitrary smooth function $f$ in terms of its Legendre transform.

\end{document}